\documentclass{amsart}
\usepackage{palatino,hhline, bbm}
\usepackage{arydshln} 
\usepackage{amsthm, amsmath}
\usepackage{amssymb} 
\usepackage{amscd}
\usepackage{mathrsfs}
\renewcommand{\mathcal}{\mathscr}
\usepackage[hmargin=3cm, vmargin=3cm]{geometry}
\usepackage[breaklinks,bookmarksopen,bookmarksnumbered]{hyperref}
\usepackage[all]{xy}

\theoremstyle{plain}
\newtheorem*{thm}{Theorem}

\newtheorem*{prop}{Proposition}
\newtheorem*{cor}{Corollary}

\theoremstyle{remark}

\newcommand\pr{\noindent\textit{Proof} : }

\newcommand\rond{\kern 1pt{\scriptstyle\circ}\kern 1pt}

\newcommand\im{\operatorname{Im}}

\newcommand\Pic{\operatorname{Pic}}

\newcommand\CH{\operatorname{CH}}

\newcommand\Z{\mathbb{Z}}

\newcommand\C{\mathbb{C}}
\renewcommand\P{\mathbb{P}}
\newcommand\N{\mathbb{N}}

\renewcommand\O{\mathcal{O}}

\newcommand\iso{\vbox{\hbox to .8cm{\hfill{$\scriptstyle\sim$}\hfill}
\nointerlineskip\hbox to .8cm{{\hfill$\longrightarrow $\hfill}} }}
\newcommand\bir{\vbox{\hbox to .8cm{\hfill{$\scriptstyle\sim$}\hfill}
\nointerlineskip\hbox to .8cm{{\hfill$\dasharrow $\hfill}} }}
\newcommand\abs[1]{\lvert {#1}\rvert}

\begin{document}
\title{A remark on the generalized Franchetta conjecture for K3 surfaces}
\author[Arnaud Beauville]{Arnaud Beauville}
\address{Universit\'e C\^ote d'Azur\\
CNRS -- Laboratoire J.-A. Dieudonn\'e\\
%UMR 7351 du CNRS\\
Parc Valrose\\
F-06108 Nice cedex 2, France}
\email{arnaud.beauville@unice.fr}
 
\begin{abstract}
A family of K3 surfaces $\mathscr{X}\rightarrow B$ has the \emph{Franchetta property} if the Chow group of 0-cycles on the generic fiber is cyclic. The generalized Franchetta conjecture proposed by O'Grady asserts that the universal family  $\mathscr{X}_g\rightarrow \mathscr{F}_g$ of polarized K3 of degree $2g-2$ has the Franchetta property. While this is known only for small $g$ thanks to \cite{PSY}, we prove that for all $g$ there is a hypersurface in $ \mathscr{F}_g$ such that the corresponding family has the Franchetta property.
\end{abstract}
%\subjclass[2010]{Primary: 14M20; Secondary: 14E08, 14J45}
\maketitle 
\section{Introduction}
In 1954, Franchetta stated that the only line bundles defined on the generic curve of genus $g\geq 2$ are the powers of the canonical bundle \cite{F}. Since the proof was insufficient, the result became known as the \emph{Franchetta conjecture}; it was proved by Harer in \cite{H}, see also \cite{AC}.

In \cite{OG}, O'Grady proposed an analogue of this result for 0-cycles on K3 surfaces. Recall that the Chow group $\CH^2(X)$ of 0-cycles on a K3 surface $X$ contains a canonical class $\mathfrak{o}_X$, the class of any point lying on some rational curve in $X$; for any divisors $D$ and $D'$ on $X$, the product $D\cdot D'$ in $\CH^2(X)$ is a multiple of $\mathfrak{o}_{X}$  \cite{BV}. Let $p:\mathscr{X}\rightarrow B$ be a map of smooth varieties whose general fiber is a K3 surface. We say that the family $\mathscr{X}\rightarrow B$ has the Franchetta property if for every smooth fiber $X$ of $p$ the image of the restriction map $\CH^2(\mathscr{X})\rightarrow \CH^2(X)$ is contained in $\Z\cdot \mathfrak{o}_X$. Equivalently, the Chow group $\CH^2(\mathscr{X}_\eta )$ of the generic fiber is cyclic.

For $g\geq 2$, let $\mathscr{X}_g\rightarrow \mathscr{F}_g$ be the universal family of polarized K3 surfaces of degree $2g-2$.
The generalized Franchetta conjecture of O'Grady is the assertion that this family 
has the Franchetta property\footnote{Here one can view $\mathscr{F}_g$ as a stack, or restrict to the open subset parametrizing K3 with trivial automorphism group.}. It is proved for $g\leq 10$ and some higher values of $g$ in \cite{PSY}; the general case seems far out of reach. We prove in this note a much weaker (and much easier) statement:
\begin{thm}
There exists for every $g$ a hypersurface in $\mathscr{F}_g$ such that the corresponding family satisfies the Franchetta property.
\end{thm}
The key point of the proof is the construction, for each $g$, of a $18$-dimensional family of polarized K3 surfaces of degree $2g-2$, which can be realized as complete intersections in $\P^1\times \P^n$ for $n=2,3$ or $4$ (\S 3). Then a simple argument, already used in \cite{PSY}, shows that these families have the Franchetta property (\S 2). Here the crucial property of our families is that they are parameterized by a linear space; thus there is no chance that the method extends to the whole moduli space $\mathscr{F}_g$, which is of general type for $g$ large enough \cite{GHS}.

\bigskip	
\section{The method}
We use the  method of \cite{PSY}, based on the following result.
Let $P$ be a smooth complex projective variety, $E$ a vector bundle on $P$, globally generated by a subspace $V$ of $H^0(E)$. Consider the subvariety $\mathscr{X}\subset \P(V)\times  P$ of pairs $(\C s,x)$ with $s(x)=0$\ \footnote{Here $\P(V)$ is the space of lines in $V$.}; let $p,q$ be the projections onto $\P(V)$ and $P$. For $s\in V\smallsetminus\{0\} $, the fiber $p^{-1}(\C s)$ is the zero locus  of $s$ in $P$; for $x\in P$, the fiber $q^{-1}(x)$ is the space of lines $\C s\subset V$ such that $s(x)=0$. Since $V$ generates $E$, the projection $q: \mathscr{X}\rightarrow P$ is a projective bundle (in particular, $\mathscr{X}$ is smooth).

\begin{prop}\label{key}
For any smooth fiber $X$ of $p$,  the image of the restriction map $\CH(\mathscr{X})\rightarrow \CH(X)$ is equal to the image of $\CH(P)$.
\end{prop}
\pr Let $h\in \CH^1(\P(V))$ be the class of a hyperplane section. The class $p^*h\in CH^1(\mathscr{X})$ induces the hyperplane class on a general fiber of $q$; since $q$ is a projective bundle, it follows that $\CH(\mathscr{X})$ is generated by $q^*\CH(P)$ and the powers of $p^*h$. But $p^*h$ vanishes on the fibers, hence the result.\qed

\begin{cor}
Assume that the smooth fibers of $p$ are K3 surfaces, and that the multiplication map 

 \noindent $m_P:\operatorname{Sym}^2  \CH^1(P)\rightarrow \CH^2(P)$ is surjective. Then the family $\mathscr{X}\rightarrow \P(V)$ has the Franchetta property.
\end{cor}
\pr Let $X$ be a smooth fiber of $p$. The commutative diagram
\[\xymatrix{\operatorname{Sym}^2  \CH^1(P) \ar[r]\ar@{->>}[d]_{m^{}_P} & \operatorname{Sym}^2  \CH^1(X)\ar[d]^{m^{}_X}\\ \CH^2(P)\ar[r]&\CH^2(X)}\]
shows that the image of $\CH^2(P)\rightarrow \CH^2(X)$ is contained in the image of $m_X$, hence in $\Z\cdot \mathfrak{o}_X$.\qed

\bigskip	
\section{Proof of the theorem}
Since $\dim \mathscr{F}_g=19$, we must construct for every $g$ a $18$-dimensional family of polarized K3 surfaces $(S,L)$ with $(L)^2=2g-2$ satisfying the Franchetta property.
 We will need three different constructions in order to cover every  $g\geq 8$  (the small genus case follow from \cite{PSY}). 
 We will apply the Corollary with $P=\P^1\times \P^n $ for $n=2,3$ or $4$ --- note that the surjectivity of $m_P$ is trivially satisfied. For $i,j\in\N$, we put $\O_P(i,j):= \O_{\P^1}(i)\boxtimes \O_{\P^n}(j)$; the vector bundle $E$ will be a direct sum of $n-1$ line bundles of this type, so $S$ is a complete intersection of $n-1$ hypersurfaces in $P$. In order for $S$ to be a K3 surface we must have $\det(E)=K_{P}^{-1}=\O_P(2,n+1)$.  We will always take $V=H^0(E)$.

The polarization $L$ on our K3 surface $S$ will be the restriction of the very ample line bundle $\O_P(a,1)$ on $P$,  for $a\geq 1$.  Let $p,h\in \CH^1(P)$ be the pull back of the class of a point in $\P^1$ and of the hyperplane class in $\P^n $. Then
\[ 2g-2=(L)^2=(ap+h)^2\cdot [S]=\bigl(2a(p\cdot h)+h^2\bigr)\cdot [S]\, .\]
\textbf{Case I :}  $n=2$,  $E=\O_P(2,3)$, hence 
\[2g-2=\bigl(2a(p\cdot h)+h^2\bigr)\cdot (2p+3h)= 2(3a+ 1)\,.\]
\textbf{Case II :}   $n=3$,   $E=\O_P(1,1)\oplus \O_P(1,3)$, hence
\[2g-2=\bigl(2a(p\cdot h)+h^2\bigr)\cdot (p+h)(p+3h)=2(3a+2)\,.\]
\textbf{Case III :}  $n=4$,   $E=\O_{P}(0,3)\oplus \O_P(1,1)\oplus \O_P(1,1)$, hence
\[2g-2=\bigl(2a(p\cdot h)+h^2\bigr)\cdot 3h(p+h)^2=2(3a+ 3)\,.\]
Thus we get all values of $g\geq 8$. 

\medskip	
 It remains to prove that the three families just constructed are $18$-dimensional. The exact sequence
\[0\rightarrow T_S\rightarrow T_{P|S}\rightarrow N_{S/P}\rightarrow 0\]
gives rise to an exact sequence \[0\rightarrow H^0(T_{P|S})\rightarrow H^0(N_{S/P})\xrightarrow{\ \partial \ } H^1(S,T_S)\,;\]
the image of $\partial $ describes, inside the space of first order deformations of $S$, those which come from our family. Thus we want to prove $\dim \im \partial =18$, or equivalently $h^0(N_{S/P})-h^0(T_{P|S})=18$. 

We have $T_P=\operatorname{pr}_1^*T_{\P^1}\oplus\, \operatorname{pr}_2^*T_{\P^n}  $; from the Euler exact sequence we get $h^0((\operatorname{pr}_1^*T_{\P^1})^{}_{|S})=h^0(\operatorname{pr}_1^*T_{\P^1})$,
 and similarly for $\operatorname{pr}_2^*T_{\P^n}$. Thus $h^0(T_{P|S})=h^0(T_{\P^1})+h^0(T_{\P^n})=3+n(n+2)$.

\smallskip	
Let us denote by $d_S$ the restriction to $S$ of a class $d\in\Pic(P)$. 
Using $d_{S}\cdot d'_{S}=d\cdot d'\cdot [S]$, we find 
\[p_S^2=0\ ,\quad p_S.h_S=3\ ,\quad h_S^2=2n-2\ .\]
By Riemann-Roch, we have $h^0(\O_S(i,j))=2+\dfrac{1}{2}(ip_S+jh_S)^2=2+3ij+ j^2(n-1)$.

\smallskip	
\textbf{Case I :} $h^0(N_{S/P})=h^0(\O_S(2,3))=29$, $h^0(T_{P|S})=11 $.

\smallskip	
 \textbf{Case II :} $h^0(N_{S/P})=h^0(\O_S(1,1))+h^0(\O_S(1,3))=9+29=36$, $h^0(T_{P|S})=18$.
 
 \smallskip	
 \textbf{Case III :} $h^0(N_{S/P})=2h^0(\O_S(1,1))+h^0(\O_S(0,3))=2\cdot 8+29=45$, $h^0(T_{P|S})=27$.
 
 \smallskip	
 In each case we find $h^0(N_{S/P})-h^0(T_{P|S}=18$ as required.\qed

\bigskip	

\noindent\emph{Remarks}$.-$ 1) In fact, for $S$ very general in each family, $\Pic(S)$ is generated by $p_S$ and $h_S$: this follows from the Noether-Lefschetz theory, see \cite[Thm. 3.33]{V}. Therefore $\Pic(S)$ is the rank 2 lattice with intersection matrix $\begin{pmatrix}
0 & 3\\ 3 & 2n-2
\end{pmatrix}$.

\medskip	
2) Our 3 families admit actually a simple geometric description. In what follows we consider a general surface $S$ in each   family. We fix homogeneous coordinates $U,V$ on $\P^1$.

\smallskip	
\textbf{Case I :}  $S$ is given by an equation $U^2A+2UV B+V^2C=0$ in $P=\P^1\times \P^2$, with $A,B,C$ cubic forms on $\P^2$.  Projecting onto $\P^2$ gives a double covering $S\rightarrow \P^2$ branched along the sextic plane curve $\Gamma :B^2-AC=0$. Let $\alpha $ and $\gamma $ be the divisors on $\Gamma $ defined by $A=B=0$ and $C=B=0$; then $2\alpha $, $2\gamma $ and $\alpha +\gamma  $ are induced by the cubic curves $A=0$, $C=0$ and $B=0$ respectively, hence belong to the canonical system $\abs{K_{\Gamma }}$. It follows that $\alpha $ and $\gamma  $ are linearly equivalent theta-characteristics, hence belong to a half-canonical $g^1_9$, that is, a vanishing thetanull on $\Gamma $. Conversely, it is easy to see that a smooth plane sextic with a  vanishing thetanull  has an equation of the above form. We conclude that \emph{the surfaces in Case I are the double covers of $\P^2$ branched along a sextic curve with a vanishing thetanull}.

\medskip	
\textbf{Case II :}  The equations of $S$ in $P=\P^1\times \P^3$ have the form $UL+VM=UA+VB=0$, where $L,M;A,B$ are forms of degree $1$ and $3$ on $\P^3$. The projection $S\rightarrow \P^3$ is an isomorphism onto the quartic surface $LB-MA=0$; this is the equation of a general quartic containing a line. Thus \emph{the surfaces in Case II are the quartic surfaces containing a line}.

\medskip	
\textbf{Case III :}  The equations of $S$ in $P=\P^1\times \P^4$ are of the form $UA+VB=UC+VD=F=0$, where $A,B,C,D;F$ are forms of degree $1$ and $3$ on $\P^3$.  The projection $S\rightarrow \P^4$ is an isomorphism onto the surface $AD-BC=F=0$, that is, the intersection of a quadric cone (with one singular point) and a cubic in $\P^4$.  Thus \emph{the surfaces in Case III are the complete intersections of a quadric cone and a cubic in $\P^4$}.

\smallskip	
Note that one sees easily from this description that each family depends indeed on 18 moduli.

\medskip	
3) Our construction shows  in particular that \emph{the moduli space $\mathscr{F}_g$ contains a unirational hypersurface for any $g$}, despite the fact that it is of general type for $g\gg 0$. 

\end{document}